\documentclass[12pt]{amsart}
\usepackage{amsmath,amscd,amssymb,amsfonts}
\newcommand{\nin}{\par\noindent}
\newcommand{\bs}{\par\bigskip}
\newcommand{\ms}{\par\medskip}
\newcommand{\sk}{\par\smallskip}
\newcommand{\msum}{\hbox{$\sum$}}
\newcommand{\mtim}{\h{$\times$}}
\newcommand{\motim}{\h{$\otimes$}}
\newcommand{\Alb}{{\rm Alb}}
\newcommand{\al}{\alpha}
\newcommand{\bl}{\bigl}
\newcommand{\br}{\bigr}
\newcommand{\C}{{\mathbf C}}
\newcommand{\dd}{\partial}
\newcommand{\fr}{{\rm fr}}
\newcommand{\ga}{\gamma}
\newcommand{\h}{\hbox}
\newcommand{\Hom}{{\rm Hom}}
\newcommand{\Ker}{{\rm Ker}}
\newcommand{\la}{\lambda}
\newcommand{\LL}{{\mathcal L}}
\newcommand{\M}{{\mathcal M}}

\newcommand{\OO}{{\mathcal O}}

\newcommand{\Pic}{{\rm Pic}}
\newcommand{\q}{\quad}
\newcommand{\R}{{\mathbf R}}
\newcommand{\om}{\omega}
\newcommand{\theo}{{}\,\overline{\!\theta}{}}
\newcommand{\X}{\widetilde{X}}
\newcommand{\Z}{{\mathbf Z}}
\newcommand{\zo}{\overline{z}}
\newcommand{\into}{\hookrightarrow}
\newcommand{\simto}{\buildrel\sim\over\longrightarrow}
\newcommand{\ssb}{\raise.25ex\h{${\scriptscriptstyle\bullet}$}}
\newcommand{\ssc}{\,\raise.15ex\h{${\scriptstyle\circ}$}\,}
\begin{document}
\title[Twistor deformation]{Twistor deformation of rank one local systems}
\author[M. Saito]{Morihiko Saito}
\address{RIMS Kyoto University, Kyoto 606-8502 Japan}
\begin{abstract}
We determine the twistor deformation of rank one local systems on compact Kaehler manifolds which correspond to smooth twistor modules of rank one in the sense of C. Sabbah. Our proof is rather elementary, and uses a natural description of the moduli space of rank one local systems together with the canonical morphism to the Picard variety.
The corresponding assertion for smooth twistor modules of rank one follows from the theory of C. Simpson and has been known to the specialists, according to T. Mochizuki.
\end{abstract}
\maketitle
\centerline{\bf Introduction}
\bs\nin
Let $X$ be a compact K\"ahler manifold of dimension $d$. Set $n:=\dim H^1(X,\OO_X)$. There is a canonical morphism
$$\rho:\Hom\bl(H_1(X,\Z)_{\fr},\C^*\br)\cong(\C^*)^{2n}\to\Pic^0(X),
\leqno(0.1)$$
where the first isomorphism is given by choosing free generators $\{\ga_i\}$ of
$$H_1(X,\Z)_{\fr}:=H_1(X,\Z)/H_1(X,\Z)_{\rm tor}.$$
Note that the source of $\rho$ is the identity component of the moduli space of rank 1 local systems
$$H^1(X,\C^*)=\Hom\bl(H_1(X,\Z),\C^*\br),$$
and the latter is the product of this identity component with a finite group by choosing a splitting of the torsion subgroup $H_1(X,\Z)_{\rm tor}$. This finite group is identified with $H^2(X,\Z)_{\rm tor}$, and can be neglected for the study of the kernel of the canonical morphism to the moduli space of line bundles.
The $\C$-local systems of rank 1 on $X$ are uniquely determined up to an isomorphism by the monodromies $\la_i\in\C^*$ along the free generators $\ga_i$, where the monodromies along torsion elements of $H_1(X,\Z)$ are assumed trivial.
\sk
The target of $\rho$ is defined by
$$\Pic^0(X):=H^1(X,\C)/\bl(\Gamma(X,\Omega_X^1)+H^1(X,\Z(1))\br).
\leqno(0.2)$$
This is the moduli space of topologically trivial holomorphic line bundles on $X$ by using the exponential sequence as is well-known, where $\Z(1):=2\pi\sqrt{-1}\,\Z$, see [De2]. (Note that $H^1(X,\Z)$ is torsion-free.) It is easy to see that
$$\h{The morphism $\rho$ is a surjective morphism of complex analytic Lie groups.}
\leqno(0.3)$$
(See Propositions~(1.1-2) below.)
\sk
Consider the morphism
$$\sigma:\Gamma(X,\Omega_X^1)\ni\theta\mapsto\bl(\exp\bl(-\h{$\int$}_{\ga_1}\theta\br),\dots,\exp\bl(-\h{$\int$}_{\ga_{2n}}\theta\br)\br)\in(\C^*)^{2n}.
\leqno(0.4)$$
In this paper, we show the following.
\ms\nin
{\bf Theorem~1.} {\it We have}
$${\rm Ker}\,\rho={\rm Im}\,\sigma.
\leqno(0.5)$$
\ms
Actually the proof of Theorem~1 turns out almost trivial (up to the verification of the compatibility of some canonical morphisms) although the minus sign in the definition of $\sigma$ in (0.4) does not appear here, see (2.2).
The first proof of Theorem~1 used the reduction to the abelian variety case via the Albanese map, and solved a certain differential equation, where we get the minus sign as in (0.4), see (2.1).
By Theorem~1 and (0.3), we get the following.
\ms\nin
{\bf Corollary~1.} {\it For any $\xi\in\Pic^0(X)$, let $\lambda=(\lambda_i)$ be an element of $\rho^{-1}(\xi)\subset(\C^*)^{2n}$. Then}
$$\rho^{-1}(\xi)=({\rm Im}\,\sigma)\,\lambda.
\leqno(0.6)$$
\ms
It is rather easy to show that there is a unique rank 1 unitary local system $\eta$ in $\rho^{-1}(\xi)$ (see Proposition~(1.3) below), and we may assume that $\lambda$ is this element.
Corollary~1 is closely related to the theory of C.~Simpson, see Example after Prop.~1.5 in [Si1], where he assumes $\dim X=1$.
His $\C^*$-action on Higgs line bundles $(\LL,\theta)$ seems to be defined by $(\LL,t\theta)$ ($t\in\C^*$) where $\theta\in\Gamma(X,\Omega_X^1)$ is the Higgs 1-form. In terms of the moduli space of rank 1 local systems, it is then given by
$$\h{$\eta\,\bl(\exp\bl(-\int_{\ga_1}(t\theta+\overline{t}\theo)\br),\dots,\exp\bl(-\int_{\ga_{2n}}(t\theta+\overline{t}\theo)\br)\br)\q(t\in\C^*),$}
\leqno(0.7)$$
where $\eta\in(\C_1^*)^{2n}$ is the rank 1 unitary local system corresponding to the line bundle $\LL$.
Note that ${\rm Re}(\theta):=(\theta+\theo)/2$ is the real part of $\theta$, i.e., the image of $\theta$ by the canonical morphism
$$\Gamma(X,\Omega_X^1)\to H^1(X,\C)=H^1(X,\R)\oplus H^1(X,\R(1))\buildrel{pr}\over\to H^1(X,\R),$$
and the exponential map is equivalent to the passage to the quotient by $H^1(X,\Z(1))$.
\sk
On the other hand, we have the twistor deformation of a rank 1 local system, i.e., a holomorphic family of rank 1 local systems on $X$ parametrized by $z\in\C^*$, which corresponds to a smooth twistor module of rank 1 in the sense of C. Sabbah [Sa] (where the local system at $z=1$ is the given local system), see also [Si2], [Si3]. The following has been informed from T.~Mochizuki [Mo]:
\ms\nin
{\bf Theorem~2.} {\it The smooth twistor module of rank $1$ associated with a Higgs line bundle $(\LL,\theta)$ can be expressed by using the holomorphic family of connections
$$\nabla_h+z^{-1}\,\theta+z\,\theo\q(z\in\C^*),$$
where $\nabla_h$ is the unitary connection on the underlying $C^{\infty}$ line bundle of $\LL$ associated with a pluri-harmonic metric $h$.}
\ms
This follows from the theory of Simpson ([Si2], Thm.~4.3) and has been known to the specialists, according to Mochizuki.
By translating Theorem~2 to an assertion on local systems (see Remark~(1.5) below), we get the following.
\ms\nin
{\bf Theorem~3.} {\it The twistor deformation, i.e. the family of rank $1$ local systems associated with the twistor module corresponding to a Higgs line bundle $(\LL,\theta)$ is written as
$$\h{$\eta\,\bl(\exp\bl(-\int_{\ga_1}(z^{-1}\,\theta+z\,\theo)\br),\dots,\exp\bl(-\int_{\ga_{2n}}(z^{-1}\,\theta+z\,\theo)\br)\br)\q(z\in\C^*).$}
\leqno(0.8)$$
where $\eta$ is as in $(0.7)$.}
\ms
We give a proof of Theorem~3 using the moduli space of rank 1 local systems together with the morphism to the Picard variety, see (2.3) below. This is rather elementary, and seems to be relatively easy to follow even for non-experts of harmonic bundles.
\sk
By Theorem~1, considering the images of (0.7) and (0.8) in the moduli space of line bundle (i.e. in the Picard variety) is equivalent to forgetting the holomorphic terms $t\theta$ and $z^{-1}\,\theta$. Hence these images coincide up to the complex conjugation of the parameters (i.e. by setting $z=\overline{t}$).
Note that $z^{-1}\,\theta+z\,\theo$ in (0.8) cannot be replaced with $z\,\theta+z^{-1}\,\theo$ as someone might think, since $z^{-1}\,\theo$ does not converge at $z=0$ if $\theta\ne 0$ (by considering its image in the Picard variety by using Theorem~1, see the proof of Theorem~3 in (2.3)).
\sk
We would like to thank those who are interested in this manuscript, especially Professors T. Mochizuki and C. Sabbah who informed us of the relation with the theory of Simpson ([Si1], [Si2], [Si3]) and also with that of twistor modules.
\ms
In Section 1 we review some well-known facts from the theory of rank one local systems. In Section 2 we prove the main theorems.
\bs\bs
\centerline{\bf 1. Preliminaries}
\bs\nin
In this section we review some well-known facts from the theory of rank one local systems.
\ms\nin
{\bf 1.1.~Proposition.} {\it The morphism $\rho$ in $(0.1)$ is a morphism of complex analytic Lie groups.}
\ms\nin
{\it Proof.} For a contractible Stein manifold $S$, the exponential sequence induces a long exact sequence
$$H^1(X,\Z(1))\to H^0(S,R^1(pr_2)_*\OO_{X\mtim S})\to H^1(X\mtim S,\OO_{X\mtim S}^*)\to H^2(X,\Z(1)),
\leqno(1.1.1)$$
by applying the Leray spectral sequence to $pr_2:X\mtim S\to S$. So a topologically trivial holomorphic line bundle on $X\mtim S$ gives a section of a locally free sheaf
$$R^1(pr_2)_*\OO_{X\mtim S}.$$
This implies that $\rho$ is complex analytic since the universal family of rank 1 local systems can be constructed by dividing the trivial line bundle on the product of the universal covering of $X$ with $(\C^*)^n$, and this gives a holomorphic family of line bundles. The compatibility with the multiplicative structures is trivial since these are defined by using the tensor product.
\ms\nin
{\bf 1.2.~Proposition.} {\it The morphism $\rho$ is surjective.}
\ms\nin
{\it Proof.} We have the commutative diagram
$$\begin{CD}0@>>>H^1(X,\Z(1))@>>>H^1(X,\C)@>>>H^1(X,\C^*)@>>>H^2(X,\Z(1))\\
@. @| @VV{\al}V @VV{\beta}V @|\\
0@>>>H^1(X,\Z(1))@>>>H^1(X,\OO_X)@>>>H^1(X,\OO_X^*)@>>>H^2(X,\Z(1))
\end{CD}
\leqno(1.2.1)$$
where $\al$ is surjective by Hodge theory. Moreover, by the universal coefficient theorem we have a canonical isomorphism
$${\rm Hom}(H_1(X,\Z),\C^*)=H^1(X,\C^*),
\leqno(1.2.2)$$
(since $\C^*$ is an injective $\Z$-module), and $\rho$ can be identified with a restriction of $\beta$. So the assertion follows.
\ms
The above argument can be modified to get the following
(which implies that there is a unique unitary local system of rank 1 in each fiber of $\rho$).
\ms\nin
{\bf 1.3.~Proposition.} {\it Set $\,\C^*_1:=\{\la\in\C^*\mid|\la|=1\}$. We have a canonical isomorphism}
$$H^1(X,\C_1^*)\simto H^1(X,\OO_X^*).$$
\ms\nin
{\it Proof.} This follows by replacing the first exact sequence in (1.2.1) with the long exact sequence associated with the exponential sequence
$$0\to\Z(1)\to\R(1)\to\C^*_1\to 1,
\leqno(1.3.1)$$
since Hodge theory implies the canonical isomorphism
$$H^1(X,\R(1))\simto H^1(X,\OO_X).
\leqno(1.3.2)$$
\ms\nin
{\bf 1.4.~Albanese maps.} For a compact K\"ahler manifold $X$ of dimension $d$, its Albanese variety $\Alb(X)$ is defined by
$$\Alb(X):=H^{2d-1}(X,\C)/\bl(F^dH^{2d-1}(X,\C)+ H^{2d-1}(X,\Z(d))_{\fr}\br),
\leqno(1.4.1)$$
where $H^{2d-1}(X,\Z)_{\fr}:=H^{2d-1}(X,\Z)/{\rm torsion}$. Choosing a point $x_0$ of $X$, the Albanese map $\al_X:X\to\Alb(X)$ can be defined by using the integrals
$$H^0(X,\Omega_X^1)\ni\om\mapsto\int_{x_0}^x\om\in\C\,\,\, {\rm mod}\,\,\,\langle H_1(X,\Z),\om\rangle,
\leqno(1.4.2)$$
where $H^{2d-1}(X,\C)/F^dH^{2d-1}(X,\C)$ is identified with the dual of $H^0(X,\Omega_X^1)$ and we have by Poincar\'e duality
$$H^{2d-1}(X,\Z(d))=H_1(X,\Z).$$
Setting $Y:=\Alb(X)$, this construction implies an isomorphism
$$(\al_X)_*:H_1(X,\Z)_{\fr}\simto H_1(Y,\Z),
\leqno(1.4.3)$$
and hence
$$\aligned\al_X^*&:H^1(Y,\Z)\simto H^1(X,\Z),\\
\al_X^*&:\Gamma(Y,\Omega_Y^1)\simto\Gamma(X,\Omega_X^1).\endaligned
\leqno(1.4.4)$$
\ms\nin
{\bf 1.5.~Remark.}
The translation between Theorems~2 and 3 consists of the calculation of the global monodromies of the $C^{\infty}$ connection
$$\nabla_h+z^{-1}\,\theta+z\,\theo.$$
By using the tensor product of connections, this can be reduced to the case $\LL=\OO_X$ and $\nabla_h=d$. Then the monodromies are calculated by using the pull-back of the connection by smooth paths representing generators of $H_1(X,\Z)_{\rm fr}$.
\bs\bs
\centerline{\bf 2. Proof of the main theorems}
\bs\nin
In this section we prove the main theorems.
\ms\nin
{\bf 2.1~First proof of Theorem~1.} By (1.4), the assertion is reduced to the case where $X$ is a complex torus of dimension $n$. We have the universal covering $\pi:\X\to X$, and $\X$ has a natural structure of a $\C$-vector space. Choosing a basis we have $\X=\C^n$. We also choose generators $\ga_i\,\,(i\in[1,2n])$ of $H_1(X,\Z)$ which correspond to generators $\om_i$ of $\Ker\,\pi\cong\Z^{2n}\subset\C^n$.
\sk
Then the problem is about an integrable holomorphic connections on a trivial line bundle on $X$, see [De1]. Let $v$ be a nonzero global section of a trivial line bundle over $X$. This is unique up to a nonzero constant multiple. Let $\xi_j\,\,(j\in[1,n])$ be the invariant vector field on $X$ with $\pi^*\xi_j=\dd/\dd x_j$ where the $x_j$ are the natural coordinates of $\C^n$. Then a holomorphic connection $\nabla$ is uniquely determined by $c_j\in\Gamma(X,\OO_X)=\C$ satisfying
$$\nabla_{\xi_j}v=c_jv\q(j\in[1,n]),
\leqno(2.1.1)$$
where the integrability trivially holds since $c_j\in\C$. This $c_j$ is independent of the choice of $v$. The connection $\nabla$ can be written as
$$\nabla=d+\theta\wedge,
\leqno(2.1.2)$$
where $\theta$ is the Higgs one-form with $\pi^*\theta=\sum_ic_idx_i$ in this case, see [Si1] for the general case of Higgs fields.
\ms
The monodromy of the corresponding local system is given by solutions of the differential equation on $X$
$$\xi_jf=-c_jf\q(j\in[1,n]),
\leqno(2.1.3)$$
which comes from
$$\nabla_{\xi_j}fv=(\xi_jf)v+c_jfv=0.$$
The pull-back of the differential equation to $\X=\C^n$ is given by
$$\dd f/\dd x_j=-c_jf.
\leqno(2.1.4)$$
This has nontrivial solutions on $\X=\C^n$ of the form:
$$f(x)=ae^{-\langle c,x\rangle}\q\h{with}\q
\langle c,x\rangle:=\msum_{j=1}^n\,c_jx_j\q\h{and}\q a\in\C^*.
\leqno(2.1.5)$$
Setting $a=1$ so that $f(0)=1$, the corresponding point of $(\C^*)^{2n}$ is then given by
$$(e^{-\langle c,\,\om_i\rangle})\in(\C^*)^{2n},
\leqno(2.1.6)$$
where the $\om_i$ are generators of $\Ker\,\pi\cong\Z^{2n}\subset\C^n$ corresponding to $\ga_i\in H_1(X,\Z)$. We thus get a family of rank 1 local systems $L_c$ on $X$ whose associated line bundles $\LL_c:=\OO_X\motim_{\C}L_c$ are trivial. This is parametrized by $c=(c_1,\dots,c_n)\in\C^n$. So Theorem~1 follows from (2.1.6). (This argument can be extended to the case of any topologically trivial line bundles by replacing the differential $d$ with the connection $\nabla$ associated to a unitary local system of rank 1.)
\ms\nin
{\bf 2.2~Simple proof of Theorem~1.} The morphisms $\sigma$ and $\rho$ can be identified respectively with the compositions of canonical morphisms
$$\Gamma(X,\Omega_X^1)\to H^1(X,\C)\to H^1(X,\C)/H^1(X,\Z(1)),
\leqno(2.2.1)$$
$$H^1(X,\C)/H^1(X,\Z(1))\to H^1(X,\C)/\bl(\Gamma(X,\Omega_X^1)+H^1(X,\Z(1))\br).
\leqno(2.2.2)$$
So the assertion follows.
\ms\nin
{\bf 2.3.~Proof of Theorem~3.} Associated with a smooth twistor module of rank 1 on $X$, we have a twistor deformation, i.e., a holomorphic family of rank 1 local systems on $X$ parametrized by $\C^*$.
(Note that a {\it smooth} twistor module means that it corresponds to a local system.) Set
$$V:=H^1(X,\C),\q\Gamma:=H^1(X,\Z(1)).$$
Any twistor deformation is expressed by using a multi-valued holomorphic function on $\C^*$ with values in $V$:
$$\Psi(z)=(2\pi i)^{-1}(\log z)\,\eta_0+\sum_k g_k(z)\,\theta_k+\sum_k h_k(z)\,\theo_k\,\,\,\,(z\in\C^*),
\leqno(2.3.1)$$
where $\eta_0\in\Gamma$, $\{\theta_k\}$ is a basis of $\Gamma(X,\Omega_X^1)$, and $g_k(z)$, $h_k(z)$ are holomorphic functions on $\C^*$.
In fact, taking the exponential map is equivalent to the passage to the quotient by $\Gamma$, see the remark after (0.7). Hence a twistor deformation is expressed by some multi-valued holomorphic function $\Psi(z)$ on $\C^*$ with values in $V$, and $\Psi(z)-(2\pi i)^{-1}(\log z)\,\eta_0$ is univalued for some $\eta_0\in\Gamma$ by considering its monodromy, since $\Psi(z)$ is well-defined modulo $\Gamma$. So (2.3.1) follows.
\sk
By the definition of twistor modules (see [Sa]) the underlying holomorphic family of line bundles is extended over $\C$. Hence $\eta_0=0$ and the $h_k(z)$ are holomorphic at $z=0$ by considering the image by the morphism to the Picard variety where the $g_k(z)\,\theta_k$ are neglected. (Indeed, $\Psi(z)$ mod $\Gamma(X,\Omega_X^1)+H^1(X,\Z(1))$ is holomorphic at $z=0$ if and only if $\Psi(z)$ mod $\Gamma(X,\Omega_X^1)$ is. Moreover, the latter condition is equivalent to that $\eta_0=0$ and the $h_k$ are holomorphic at $z=0$.)
As for the $g_k$, we see that the $zg_k(z)$ are holomorphic at $z=0$ by the definition of twistor modules in loc.~cit.\ (considering the corresponding holomorphic family of connections).  We thus get
$$\eta_0=0,\q g_k=\sum_{i\ge-1} g_{k,i}\,z^i,\q h_k=\sum_{i\ge 0} h_{k,i}\,z^i.$$
\sk
From the polarizability of twistor module which is a kind of Hermitian self-duality, we can deduce
$$\sum_k g_k(z)\,\theta_k+\sum_k h_k(z)\,\theo_k=-\sum_k \overline{g_k(-1/\zo)}\,\theo_k-\sum_k \overline{h_k(-1/\zo)}\,\theta_k\q{\rm mod}\,\,\,\Gamma.\leqno(2.3.2)$$
Here the right-hand side is obtained by calculating the Hermitian dual (which produces the complex conjugation on the values together with the minus sign) and applying the ``complex conjugation" defined by $z\mapsto-1/\zo$, see loc.~cit.
\sk
Set
$$\zeta:=\sum_k g_{k,0}\,\theta_k,\q\zeta':=\sum_k h_{k,0}\,\theo_k,\q\xi:=\sum_k g_{k,-1}\,\theta_k,\q\xi':=\sum_k h_{k,1}\,\theo_k.$$
Then (2.3.2) implies
$$\overline{g}_{k,-1}=h_{k,1},\q g_{k,i}=0\,\,\,(i\ge 1),\q h_{k,i}=0\,\,\,(i\ge 2),$$
and
$$\zeta'=-\overline{\zeta}\,\,\,{\rm mod}\,\,\,\Gamma,\q\xi'=\overline{\xi},$$
since ${\rm mod}\,\,\Gamma$ affects only the constant term.
We thus get
$$\sum_k g_k(z)\,\theta_k+\sum_k h_k(z)\,\theo_k=\zeta-\overline{\zeta}+z^{-1}\,\xi+z\,\overline{\xi}\q{\rm mod}\,\,\,\Gamma.
\leqno(2.3.3)$$
Here $\zeta-\overline{\zeta}\,\,({\rm mod}\,\,\,\Gamma)$ is identified with a rank 1 unitary local system $\eta$. This corresponds to the line bundle $\LL$ by considering the limit for $z\to 0$ in the Picard group (by using Theorem~1) and comparing this with (0.7). Moreover $\xi$ coincides with the Higgs 1-form $\theta$ of the Higgs line bundle $(\LL,\theta)$ by comparing (0.7) and (0.8) (with $\theta$ replaced by $\xi$) at $z=1$. So the assertion follows from (2.3.3).
\ms\nin
{\bf 2.4.~Hodge spectral sequence.}~For a local system $L_c$ in (2.1) (which is identified with a local system on the given K\"ahler manifold $X$ by (1.4.3)), we have the Hodge spectral sequence
$$E_1^{p,q}=H^q(X,\Omega_X^p\motim_{\OO_X}\LL_c)\Rightarrow H^{p+q}(X,L_c),
\leqno(2.4.1)$$
where $\LL_c:=\OO_X\otimes_{\C}L_c$.
This does not degenerate at $E_1$, for instance, if $c$ is nonzero and belongs to the kernel of $\rho$ in the notation of Theorem~1, since $H^j(X,L_c)=0$ for any $j$. Indeed, $H^j(X,L_c)$ is calculated by using the Koszul complex associated with the multiplication by $\la_i-1\,\,(i\in[1,2n])$ on $\C$ where $\la_i=e^{-\langle c,\,\om_i\rangle}$ is the monodromy along the path $\ga_i$.
(For the proof of the $E_1$-non-degeneration, it is actually enough to compare $H^0(X,\LL_c)$ and $H^0(X,L_c)$ in this case.)
The $E_1$-differential $d_1$ is given by the Higgs one-form $\theta$, and we have
$$\pi^*\theta=\msum_i\,c_idx_i.$$
Hence the $E_2$-term vanishes in this case, and in particular, the spectral sequence degenerates at $E_2$.
\sk
It seems that this $E_2$-degeneration of (2.4.1) holds for any rank~1 local systems $L$ on smooth projective varieties (see [Ar1], [Ar2]) by comparing the spectral sequence (2.4.1) (with $L_c$, $\LL_c$ replaced by $L$ and $\LL:=\OO_X\otimes_{\C}L$) with the following spectral sequence:
$$E_1^{p,q}=H^q(X,\Omega_X^p\motim_{\OO_X}\LL)\Rightarrow H^{p+q}(X,(\Omega_X^{\ssb}\motim_{\OO_X}\LL,\theta\wedge)).
\leqno(2.4.2)$$
These two spectral sequences are induced by the truncations $\sigma_{\ge p}$ on $\Omega_X^{\ssb}$ (see [De2]), and are associated with the complexes
$$(\Omega_X^{\ssb}\motim_{\OO_X}\LL,\nabla+\theta\wedge),\q (\Omega_X^{\ssb}\motim_{\OO_X}\LL,\theta\wedge),
\leqno(2.4.3)$$
which have the same components but different differentials
$$\nabla':=\nabla+{\theta\wedge}\q\h{and}\q\theta\wedge,$$
where $\nabla$ is the unitary connection on $\LL:=L\motim\OO_X$. Recall that the first complex is the de Rham complex associated with $(\LL,\nabla')$, where $\nabla'=\nabla+\theta\wedge$ is a non-unitary connection on a line bundle $\LL$ (assuming $\theta\ne 0)$. This complex is quasi-isomorphic to a rank 1 local system which will be denoted by $L_{\theta}$.
\sk
The two spectral sequences have the same $E_2$-terms since the morphism induced by $\nabla$ in the $E_1$-differential $d_1$ vanishes by applying Hodge theory to the unitary connection $\nabla$. Moreover, the second spectral sequence for the differential $\theta\wedge$ degenerates at $E_2$, see Prop.~3.7 and Remark~1 after its proof in [GL]. So the problem is whether they have the same total dimensions of the $E_{\infty}$-terms. (Here Lemma~2.2 in [Si1] does not immediately imply this by the same reason as in Remark~(2.5)(i) below.) Note that the $E_2$-degeneration implies that $H^j(X,L_{\al\theta})$ is independent of $\al\in\C^*$ since this holds for the complex defined by $\al\theta\wedge$.
\ms\nin
{\bf 2.5.~Remarks.}~(i) The connection $\nabla=d+\theta\wedge$ in (2.1.2) associated with a Higgs one-form $\theta$ is a {\it holomorphic connection} on a trivial holomorphic line bundle, and the {\it $C^{\infty}$ connection} on the trivial $C^{\infty}$ line bundle which is associated to $\theta$ in [Si1] does not necessarily have the same solution local system as that of the above holomorphic connection $\nabla$. This is informed from Professor T.~Mochizuki, and we would like to thank him.
As is seen from Example after Prop.~1.5 in loc.~cit., the {\it $C^{\infty}$ connection} on the trivial line bundle which is associated to a Higgs one-form $\theta$ in loc.~cit.\ seems to be $d+(\theta+\bar{\theta})\wedge$, where the latter $d$ is a {\it $C^{\infty}$ connection}, and is the sum of $\partial$ and $\bar{\partial}$.
\ms
(ii) By Proposition~(1.3) there is a unique unitary local system of rank 1 in each fiber of $\rho$. This can be shown also by using the fact that the local systems $L_c$ in (2.1) are never unitary local systems, i.e.\ they do not belong to $(\C^*_1)^{2n}$. (Indeed, $e^{-\langle c,\,\om_i\rangle}$ belongs to $\C^*_1$ if and only if the real part of $\langle c,\,\om_i\rangle$ vanishes.) This uniqueness implies that we have a section of $\rho$, but this is never complex analytic.
\ms
(iii) An argument similar to the proof of Theorem~3 in (2.3) seems to apply to the case $X=\C^*$, where the local systems of rank 1 on $X$ are parametrized by $\C^*$ via the monodromy eigenvalue. The twistor deformation in this case seems to be expressed as
$$z^m\exp(i\alpha+z\beta+z^{-1}\overline{\beta})\,\,\,(z\in\C^*)\q\h{for}\,\,\,m\in 2\Z,\,\,\alpha\in\R,\,\,\beta\in\C.$$
Here it is unclear whether $m=0$ and $\beta$ is a real number unless we assume that the twistor module is extendable over a partial compactification of $\C^*$ so that the nearby cycle functor at this point can be considered. (Indeed, the condition on the roots of the $b$-function of twistor deformation implies that $m=0$ and $\beta$ is a real number.)

\ms
(iv) In twistor theory, it seems that the twistor module $\M_X$ corresponding to the constant sheaf on a smooth complex manifold $X$ has always weight 0 independently of the dimension of $X$, see [Sa]. However, this does not seem to be compatible, for instance, with the calculation of the nearby cycles in the normal crossing case in loc.~cit., Lemma~3.7.9. In fact, it seems to assert that the primitive part of the $\ell$-th graded piece of the monodromy filtration of the nearby cycles of $\M_X$ for the function $f=x_1x_2$ is given by
$$P{\rm Gr}^W_{\ell}\psi_f\M_X=\begin{cases}\bigoplus_{a=1,2}\,(i_{Y_a})_+\M_{Y_a}&\h{if}\,\,\,\,\ell=0,\\(i_Z)_+\M_Z&\h{if}\,\,\,\,\ell=1,\end{cases}$$
where $Y_a:=\{x_a=0\}$, $Z=Y_1\cap Y_2$. (Here the monodromy is unipotent, and $\psi_{f,-1}$ is denoted by $\psi_f$ to simplify the notation.)
Indeed, we get the above formula from the one in loc.~cit., if we set
$$\M_{Y_a}=\M_X/x_a\M_X,\quad\M_Z=\M_X/(x_1,x_2)\M_X,$$
see Lemma~3.7.8 in loc.~cit., where $x_1x_2\M$ seems to mean rather $(x_1,x_2)\M$.
However, it is quite unclear why we get {\it naturally} the shift of weight for $\ell=1$. Notice that the weight of a twistor module can be changed {\it arbitrarily} since the Tate twist $(i/2)$ for $i\in\Z$ exists in twistor theory. However, the Tate twist $(-1/2)$ should be noted explicitly after $(i_Z)_+\M_Z$ in the above formula for $\ell=1$, and it should be clarified where this Tate twist comes from.
\sk
There is also a problem about the shift of the monodromy filtration on the nearby and vanishing cycle functors. In the minimal extension case, the vanishing cycle functor $\varphi_f$ is identified with the image of $N$ on the nearby cycle functor $\psi_f$ (restricted to the unipotent monodromy part). Then the monodromy filtration on the nearby cycles must be shifted by 1 in order that the morphism ${\rm can}:\psi_f\to\varphi_f$ induce isomorphisms between the primitive parts of the graded pieces of the monodromy filtration (except for the lowest degree). In fact, this is clearly impossible without shifting the center of the symmetry for one of them since the canonical morphism ``can'' must {\it preserve the weights}.
But the filtration on the vanishing cycles cannot be shifted by considering the case where the module is supported inside the divisor.
So we have to shift the monodromy filtration on the nearby cycle functor by 1.
\sk
The above problems can be solved rather naturally if one thinks, for instance, that the weights decrease by the codimension under the direct images by closed immersions. In the above case it is natural to put the Tate twist $(-1)$ after $\M_Z$, since we get $\M_Z$ for the co-primitive part naturally. This is closely related with the inclusion of twistor modules
$$\M_Y\into\psi_f\M_X,$$
which appears in the exact sequence
$$\aligned0\to\M_Y\to\psi_f\M_X\to\varphi_f\M\to 0,\\
\h{or}\q0\to\M_Y\to\psi_f\M_X\buildrel{N}\over\to\psi_f\M(-1),\endaligned$$
where $\M_Y$ is the twistor module on $X$ corresponding to the constant sheaf on $Y:=Y_1\cup Y_2$. (It may coincide with $i_+i^*\M_X$ up to a shift of complex if $i^*$ can be defined appropriately, where $i:Y\into X$ is the inclusion.)
\sk
It is expected that the weight filtration $W$ on $\M_Y$ satisfies:
$${\rm Gr}^W_k\M_Y=\begin{cases}\bigoplus_{a=1,2}\,(i_{Y_a})_+\M_{Y_a}&\h{if}\,\,\,\,k=-1,\\(i_Z)_+\M_Z&\h{if}\,\,\,\,k=-2.\end{cases}$$
In fact, this should be closely related with
$${\rm Gr}^W_kj_!\M_U=\begin{cases}\M_X&\h{if}\,\,\,\,k=0,\\
\bigoplus_{a=1,2}\,(i_{Y_a})_+\M_{Y_a}&\h{if}\,\,\,\,k=-1,\\(i_Z)_+\M_Z&\h{if}\,\,\,\,k=-2,\end{cases}$$
where $U:=X\setminus Y$ with the inclusion $j:U\into X$. Indeed, there should be a short exact sequence of twistor modules on $X$
$$0\to\M_Y\to j_!\M_U\to\M_X\to 0.$$
Note, however, that there is no nontrivial morphism of twistor modules on $X$:
$$\M_X\to\M_Y,$$
and we have only an element in the extension class, i.e.
$$\M_X\to\M_Y[1].$$
This comes from the difference in $t$-structure for constructible sheaves and $D$-modules.

\end{document}